\documentstyle[12pt,amsfonts,amssymb,leqno]{article}

\newfont\got{eufm10}

\newtheorem{proposition}{Proposition}[section]
\newtheorem{thm}[proposition]{Theorem}

\newtheorem{lemma}[proposition]{Lemma}
\newtheorem{defn}[proposition]{Definition}

\newcounter{secnum}

\begin{document}
\setcounter{section}{-1}

\begin{center}
{\Large \bf A universal weasel} 
\end{center}
\begin{center}
{\Large \bf without large cardinals in $V$} 
\end{center}
\begin{center}
\renewcommand{\thefootnote}{\fnsymbol{footnote}}
\renewcommand{\thefootnote}{arabic{footnote}}
\renewcommand{\thefootnote}{\fnsymbol{footnote}}
{\large Ralf-Dieter Schindler}
\renewcommand{\thefootnote}{arabic{footnote}}
\end{center}
\begin{center} 
{\footnotesize
{\it Institut f\"ur formale Logik, Universit\"at Wien, 1090 Wien, Austria}} 
\end{center}

\begin{center}
{\tt rds@logic.univie.ac.at}

{\tt http://www.logic.univie.ac.at/${}^\sim$rds/}\\
\end{center}

\section{Introduction.}

In \cite{CMIP}, Steel constructs an $\Omega+1$ iterable premouse, called $K^c$, 
of height
$\Omega$ which is universal in
the sense that it wins the coiteration against every coiterable premouse of height
$\leq \Omega$. Here, $\Omega$ is a fixed measurable cardinal, and Steel works in the
theory ``$ZFC + \Omega$ is measurable $+$ there's no inner model with a Woodin
cardinal.'' In \cite{NFS}, Jensen shows that ``$\Omega$ is measurable'' can
be relaxed to ``$\Omega$ is inaccessible'' here. Universal weasels are
needed for the purpose of isolating $K$, the core model.

It would be desirable to replace $\Omega$ by $OR$ here, where $OR$ is
the class of all ordinals, and to get rid of having to assume $\Omega$ ($OR$, that
is) to be ``large.'' I.e., we would like 
to prove the existence of a universal weasel in the theory 
``$ZFC +$ there's no inner model with a Woodin
cardinal.'' This would be a first step towards proving the existence of $K$ in that
theory (cf. the discussion in the introduction to \cite{FSIT}).

This note solves the problem of constructing a universal weasel. 
We prove:

\begin{thm}\label{main-thm} 
Assume $ZFC \ + $ there's no inner model with a Woodin
cardinal. There is then a universal weasel.
\end{thm}

We warn the reader that some care is necessary in order to arrive at the appropriate
notion of ``universal'' so as to make \ref{main-thm} not false for the wrong
reasons (this has to do with iteration trees of length $OR$, and will be discussed
below). 

The key new idea here is to weaken the concept of ``countably certified'' of
\cite{CMIP} Def. 1.2 which is crucial for the construction of $K^c$. 
Whereas the iterability proof of \cite{CMIP} can be checked to still go through with
this weaker requirement on new extenders to be added to the $K^c$-sequence, an
argument from \cite{NFS} can be varied to prove the universality of $K^c$. 

We do not know whether the $K^c$ constructed here satisfies a useful version of
weak covering. (It can be
shown that it does not necessarily satisfy weak covering at every (countably closed)
singular cardinal.)

\section{The existence of $K^c$.}

We build upon \cite{FSIT} and \cite{CMIP}; in particular, we use the concept of
``premouse'' as isolated there. 

Let ${\cal L} = \{ {\dot \epsilon},{\dot U} \}$ be the language for structures of the
type $(N;E,U)$, where $E$ is binary and $U$ is unary. For the purposes of this paper
let us introduce the following.

\begin{defn}
Let $\Psi$ be a (first order) formula in the language ${\cal L}$. Then $\Psi$ is
said to be $r\Sigma_3$ ({\em restricted} $\Sigma_3$) iff $$\Psi \equiv 
\exists v_0 \ ( \Phi_0 \wedge \ \Phi_1 ) {\rm , }$$
where $\Phi_0$ is $\Sigma_2$ in the language $\{ {\dot \epsilon} \}$, and 
$\Phi_1$ is $\Sigma_0$ in the language ${\cal L}$. 
\end{defn}

${\cal L}$ will be an appropriate language for the models witnessing certifiability.    

\begin{defn}\label{certificate}
Let ${\cal M}$ be an active premouse, 
$F$ the extender coded by ${\dot F}^{\cal M}$, 
$\kappa = c.p.(F)$, and $\nu = {\dot \nu}^{\cal M} =$ the natural length of $F$. 
We say that ${\cal M}$ (or, $F$) 
is {\em countably certified} iff for all ${\vec A} = (A_k
\colon k<\omega)$ with $\forall k \ \exists n = n(k) 
\ A_k \in {\cal P}([\kappa]^n) \cap
{\cal M}$ there are $N'$ and ${\vec B} = (B_k
\colon k<\omega)$ such that

(a) ${}^\omega V_\kappa \subset V_\kappa$ (i.e., $cf(\kappa) > \omega$), 
$N'$ is transitive,
${}^\omega N' \subset N'$, and $V_{\nu+1} \subset N'$, 

(b) $(V_\kappa;\in,{\vec A}) \prec_{r\Sigma_3} (N';\in,{\vec B})$, and

(c) for all $k<\omega$, $B_k \cap [\nu]^{n(k)} = i_F(A_k) \cap [\nu]^{n(k)}$, where
$i_F \colon {\cal M} \rightarrow_{\Sigma_1} Ult_0({\cal M},F)$ is the canonical
embedding.
\end{defn}

In this definition, we confuse ${\vec A}$ (and ${\vec B}$) with
$\{ (k,u) \colon u \in A_k \}$ (and $\{ (k,u) \colon u \in B_k \}$, resp.).

It is easy to see that if ${\cal M}$ is countably certified in the sense of
\cite{CMIP} Def. 1.2 then it is countably certified in the sense of
\ref{certificate}. 

We construct the models ${\cal N}_\xi$ and ${\cal M}_\xi$ as on p. 6 f. of
\cite{CMIP}, except that we don't require $(2)$ and $(3)$ at all 
in {\em Case 1} and that in
$(1)$ of {\em Case 1} we understand ``countably certified'' in the sense of
\ref{certificate} rather than \cite{CMIP} Def. 1.2.

We now have to prove \cite{CMIP} Thm. 2.5, the assertion that if ${\cal N}_\theta$
exists then collapses of countable submodels of ${\frak C}_k({\cal N}_\theta)$
are countably iterable for every $k<\omega$
(cf. \cite{CMIP} for the exact statement). As on pp. 12 ff.
we'll prove this in a simplified case, for trees of length $\omega$. We'll leave it as
an easy exercise for the reader to check that the proof of \cite{FSIT} Thm. 9. 14
can be varied in much the same way as the proof from \cite{FSIT} pp. 12 ff. in the
light of our new meaning of ``countably certified.''  
 
\begin{lemma}\label{warmup}
Let $\sigma \colon {\cal P} \rightarrow {\cal N}_\eta$ with ${\cal N}_\eta \models
ZFC$, and let ${\cal T}$ be an
iteration tree on ${\cal P}$ of length $\omega$ such that ${\cal D}^T = \emptyset$.
Then there are $b$ and $\sigma'$ such that $b$ is a cofinal branch through $T$ and
$\sigma'
\colon {\cal M}^{\cal T}_b \rightarrow {\cal N}_\eta$ with $\sigma' \circ \pi^{\cal
T}_{0b} = \sigma$.
\end{lemma}

{\sc Proof.} For any $\tau \colon {\cal P} \rightarrow {\cal Q}$ we denote by 
$$U(\tau,{\cal Q})$$ the tree of attempts to find $b$, $\tau'$ 
such that $b$ is a cofinal branch through $T$ and
$\tau'
\colon {\cal M}^{\cal T}_b \rightarrow {\cal Q}$ with $\tau' \circ \pi^{\cal
T}_{0b} = \tau$. We let 
$U(\tau,{\cal Q})$ consists of $\phi \colon {\cal M}^{\cal T}_i
\rightarrow {\cal Q}$, and if $\phi \colon {\cal M}^{\cal T}_i
\rightarrow {\cal Q}$ and $\phi' \colon {\cal M}^{\cal T}_k
\rightarrow {\cal Q}$ then we put $\phi \leq_{U(\tau,{\cal Q})} \phi'$ iff
$i \leq_T k$ and $\phi' \circ \pi^{\cal T}_{ik} = \phi$.

Let us assume that $U(\sigma,{\cal N}_\eta)$ is well-founded
(in the obvious sense). We aim to derive a
contradiction. Let us write ${\cal P}_i$ for ${\cal M}^{\cal T}_i$, and $\pi_{ik}$ for
$\pi^{\cal T}_{ik}$ (if $i \leq_T k$). Set $\kappa_i = c.p.(E^{\cal T}_i)$, and
$\nu_i =$ natural length of $E^{\cal T}_i$.

We closely follow \cite{CMIP} p. 12 ff. We are going to define 
$$(\sigma_i , {\cal Q}_i , {\cal R}_i \ \colon \ i<\omega)$$ such that the following
requirements are met, for all $j < i < \omega$. (In what follows,
the $\tau_{(\bullet,\bullet)}$'s are the functions from \cite{K} Lemma 3.1.)

\bigskip
$(1)$ ${\cal R}_i$ is a transitive model of $ZFC^-$ with ${}^\omega {\cal R}_i
\subset {\cal R}_i$,

$(2)$ $\sigma_i \colon {\cal P}_i \rightarrow {\cal Q}_i$, where ${\cal Q}_i$ is an
``${\cal N}$-model'' of ${\cal R}_i$,

$(3)$ $V^{{\cal R}_j}_{\sigma_j(\nu_j)+1} = V^{{\cal R}_i}_{\sigma_j(\nu_j)+1}$
and $\sigma_j(\nu_j) \leq \sigma_i(\nu_i)$,

$(4)$ $\tau_{(\bullet,\bullet)} \circ 
\sigma_j \upharpoonright \nu_j = \sigma_i \upharpoonright \nu_j$,

$(5)$ if $U = U(\sigma_i \circ \pi_{0i},{\cal Q}_i)$ then $U$ is well-founded and
${\cal R}_i$ has (in order type) at least $|\sigma_i|_U$ many cutoff points, and

$(6)$ $i>0 \Rightarrow {\cal R}_i \in {\cal R}_{i-1}$.

\bigskip
It is $(6)$ which gives the desired contradiction. 

To commence, we let $\sigma_0 = \sigma$, ${\cal Q}_0 = {\cal N}_\eta$, and ${\cal R}_0
= H_\theta$ for some large enough $\theta$. 

Suppose now we are given $(\sigma_l , {\cal Q}_l , {\cal R}_l \ \colon \ l \leq i)$. We
want to construct $\sigma_{i+1}$, ${\cal Q}_{i+1}$, and ${\cal R}_{i+1}$. 

Set $F = \sigma_i(E^{\cal T}_i)$, $\kappa = c.p.(F)$, and $\nu =$ the natural length
of $F$. 
Let us cheat by assuming $F$ is the top extender of ${\cal Q}_i$.
(If not, we have to consider the top extender of the target
model of $\tau_{(\bullet,\bullet)} \colon {\cal J}^{{\cal
Q}_i}_{lh(F)} \rightarrow {\tilde {\cal Q}}$ instead; a similar cheating appears in
\cite{CMIP} p. 12 ff.)
By $(2)$, $F$ is countably
certified inside ${\cal R}_i$. Let ${\vec A} = (A_k \colon k<\omega)$ be an
enumeration of $${\cal P}([\kappa]^{<\omega}) \cap {\cal Q}_i \cap ran(\sigma_i).$$
By $(1)$, ${\vec A} \in {\cal R}_i$, and hence there are inside ${\cal R}_i$ objects
$N$, $N'$, and ${\vec B}$ such that

\bigskip
(a) $N = V_\kappa^{{\cal R}_i}$, ${}^\omega N \subset N$, $N'$ is transitive,
${}^\omega N' \subset N'$, and $V_{\nu+1}^N \subset N'$, 

(b) $(N;\in,{\vec A}) \prec_{r\Sigma_3} (N';\in,{\vec B})$, and

(c) for all $k<\omega$, $B_k \cap [\nu]^{<\omega} = i_F(A_k) \cap [\nu]^{<\omega}$.

\bigskip
It now clearly suffices to prove the following 

\bigskip
{\bf Main Claim.} In $N'$, there are $\sigma$, ${\cal Q}$, and ${\cal R}$ such that

$(1)$' ${\cal R}$ is a transitive model of $ZFC^-$ with ${}^\omega {\cal R}
\subset {\cal R}$,

$(2)$' $\sigma \colon {\cal P}_{i+1} \rightarrow {\cal Q}$, where ${\cal Q}$ is an
``${\cal N}$-model'' of ${\cal R}$,

$(3)$' $V^{\cal R}_{\sigma_i(\nu_i)+1} = V^{N'}_{\sigma_i(\nu_i)+1}$ (and
hence $= V^{{\cal R}_i}_{\sigma_i(\nu_i)+1}$),

$(4)$' $\sigma \upharpoonright \nu_i = \sigma_i \upharpoonright \nu_i$, and

$(5)$' if $U = U(\sigma \circ \pi_{0i+1},{\cal Q})$ then $U$ is well-founded and
${\cal R}$ has (in order type) at least $|\sigma|_U$ many cutoff points.

\bigskip
Notice that the assertion of the Main Claim is $\Sigma_2^{N'}(\{ {\cal T},
\sigma_i(\nu_i) , \sigma_i \upharpoonright \nu_i \})$. Let $$\pi \ \colon \
(M;\in,{\vec C}) \rightarrow (N';\in,{\vec B}) {\rm , }$$ where $M$ is
countable (and hence $(M;\in,{\vec C}) \in N' \cap N$), and $\pi$ is
$\Sigma_2$-elementary w.r.t. $\{ {\dot \epsilon} \}$ and $\Sigma_0$-elementary w.r.t.
${\cal L}$. The fact that $(M;\in,{\vec C})$
can be embedded into $N'$ in such a fashion is a $r\Sigma_3$-fact, and hence
by $(N;\in,{\vec A}) \prec_{r\Sigma_3} (N';\in,{\vec B})$ there is some 
$$\pi' \ \colon \
(M;\in,{\vec C}) \rightarrow (N;\in,{\vec A})$$ 
such that $\pi'$ is
$\Sigma_2$-elementary w.r.t. $\{ {\dot \epsilon} \}$ and $\Sigma_0$-elementary w.r.t.
${\cal L}$.
In order to finish the
proof of the Main Claim (and thus of \ref{warmup}), it now suffices to verify the
following 

\bigskip
{\bf Claim.} In $N$, there are $\sigma$, ${\cal Q}$, and ${\cal R}$ such that

$(1)$'' ${\cal R}$ is a transitive model of $ZFC^-$ with ${}^\omega {\cal R}
\subset {\cal R}$,

$(2)$'' $\sigma \colon {\cal P}_{i+1} \rightarrow {\cal Q}$, where ${\cal Q}$ is an
``${\cal N}$-model'' of ${\cal R}$,

$(3)$'' $V^{\cal R}_{\pi' \circ \pi^{-1}(\sigma_i(\nu_i))+1} = 
V^{N}_{\pi' \circ \pi^{-1}(\sigma_i(\nu_i))+1}$,

$(4)$'' $\sigma \upharpoonright \nu_i = \pi' \circ \pi^{-1}(\sigma_i 
\upharpoonright \nu_i)$, and

$(5)$'' if $U = U(\sigma \circ \pi_{0i+1},{\cal Q})$ then $U$ is well-founded and
${\cal R}$ has (in order type) at least $|\sigma|_U$ many cutoff points.

\bigskip
Let $j = T$-pred$(i+1)$. We define $\sigma' \colon {\cal P}_{i+1} \rightarrow {\cal
Q}_j$ by $$\pi_{ji+1}(f)(a) \mapsto \sigma_j(f)(\pi' \circ \pi^{-1}(\sigma_i(a_i))).$$
To see that this is well-defined and elementary we argue as follows.

\begin{eqnarray*}
& {\cal P}_{i+1} \models \Phi(\pi_{ji+1}(f)(a)) \\
\Leftrightarrow & \{ u \colon {\cal P}_j \models \Phi(f(u)) \} \in (E_i^{\cal T})_a \\ 
\Leftrightarrow & \sigma_i(\{ u \colon {\cal P}_j \models \Phi(f(u)) \}) 
\in F_{\sigma_i(a)}. \\
\end{eqnarray*}
Let $\sigma_i(\{ u \colon {\cal P}_j \models \Phi(f(u)) \}) = A_k$. So $i_F(A_k) \cap
[\nu]^{<\omega} = B_k \cap
[\nu]^{<\omega}$, and we may continue as follows.
\begin{eqnarray*}
\Leftrightarrow & \sigma_i(a) \in i_F(A_k) \\
\Leftrightarrow & \sigma_i(a) \in B_k \\
\Leftrightarrow & \pi' \circ \pi^{-1}(\sigma_i(a)) \in A_k. \\ 
\end{eqnarray*}
But $A_k = \sigma_i(\{ u \colon {\cal P}_j \models \Phi(f(u)) \}) =
\sigma_j(\{ u \colon {\cal P}_j \models \Phi(f(u)) \}) =
\{ u \colon {\cal Q}_j \models \Phi(\sigma_j(f)(u)) \}$, and hence
\begin{eqnarray*}
\Leftrightarrow & {\cal Q}_j \models \Phi(\sigma_j(f)(\pi' \circ
\pi^{-1}(\pi_i(a)))). \\ 
\end{eqnarray*}

We'll have that $\sigma' \circ \pi_{0 i+1} = \sigma_j \circ \pi_{0j}$, and so
$\sigma' \in U(\sigma_j \circ \pi_{0 j},Q_j)$. Moreover, clearly,
$$\epsilon = |\sigma'|_{U(\sigma_j \circ \pi_{0 j},Q_j)} <
|\sigma_j|_{U(\sigma_j \circ \pi_{0 j},Q_j)} {\rm , }$$ and hence
by $(5)$ we may let $\theta =$ the $\epsilon^{th}$ cutoff point of ${\cal R}_j$.
Working inside ${\cal R}_j$, we may thus set

\begin{eqnarray*}
{\cal R} & = & {\rm \ the \ transitive \ collapse \ of \ the \ closure \
of \ } V_{\pi' \circ
\pi^{-1}(\sigma_i(\nu_i))+1} \cup \{ {\cal Q}_j \} \\
& & {\rm \ under \ Skolem \ functions \ for \ } V_\theta^{{\cal R}_j} 
{\rm \ and \ } \omega{\rm
-sequences,} \\
{\cal Q} & = & {\rm \ the \ image \ of \ } {\cal Q}_j {\rm \ under \ the \ collapse, \
and } \\
\sigma & = & {\rm \ the \ image \ of \ } \sigma' {\rm \ under \ the \ collapse.} \\
\end{eqnarray*}

It is now straightforward that we have shown the Claim.

\bigskip
\hfill $\square$ (\ref{warmup})

\bigskip
Of course by standard arguments the previous sketch also shows that $K^c$ exists
unless there is a non-tame premouse, say.

\section{$K^c$ is universal.}

Assume that there is no inner model with a Woodin cardinal. By the results in \S 1
together with \cite{CMIP} Lemma 2.4 (b) we then have that $K^c$ is $< OR$ iterable.
However, it may be the case that there is a definable tree on $K^c$ of length $OR$
with no cofinal branch.

This discussion leads us to the following.

\begin{defn}\label{def-universal}
A weasel $W$ is {\em universal} iff whenever $({\cal T},{\cal U})$ is a coiteration of
$W$ with some premouse ${\cal M}$ (using padded trees) with $lh({\cal T}) = lh({\cal
U}) = OR+1$ then ${\cal M}$ is a weasel, ${\cal D}^{\cal U} \cap [0,OR]_U =
\emptyset$, $\pi^{\cal U}_{0 OR}$ ''$OR \subset OR$, and ${\cal M}_{OR}^{\cal U}
\trianglerighteq {\cal M}_{OR}^{\cal T}$.
\end{defn}

N.B.: ``$W$ is universal'' is a schema which cannot be expressed by a single sentence
in the language of $ZFC$.

I do not know if there is a notion of ``universal'' which is more useful.

Let us say that a premouse ${\cal M}$ is {\em below superstrong} iff for all $F =
E^{\cal M}_\alpha \not= \emptyset$ we have that the natural length of $F$ is strictly
less than $i_F(c.p.(F))$. 
We're now going to show:

\begin{thm}\label{Kc-is-universal}
Assume $ZFC +$ every premouse is below superstrong. Then $K^c$ is universal, if it
exists.
\end{thm}

{\sc Proof.} Deny. Set $W = K^c$, and $W_\alpha = {\cal J}^W_{\alpha^{+W}}$ for
$\alpha \in OR$. By a slight refinement (due to Zeman and the author)
of an argument of Jensen (cf. \cite{NFS}) there is then a (definable)
class $C \subset OR$, club in $OR$, together with a commuting system 
$(\pi_{\alpha \beta} \ \colon \ \alpha
\leq \beta \in C)$ of maps such that
for all $\alpha
\leq \beta \in C$ do we have that
$\pi_{\alpha \beta} \colon W_\alpha \rightarrow W_\beta$ is cofinal with 
$\pi_{\alpha \beta} \upharpoonright \alpha = id$ and $\pi_{\alpha \beta}(\alpha) =
\beta$, and such that $(W_\alpha,\pi_{\alpha \beta} \ \colon \ \alpha \leq \beta \in C
\cap \kappa+1)$ is the direct limit of 
$(W_\alpha,\pi_{\alpha \beta} \ \colon \ \alpha \leq \beta \in C
\cap \kappa)$ for limit points $\kappa$ of $C$.

Let $n<\omega$ be large enough. There is then a (definable) $C' \subset C$, again club
in $OR$, such that $$N^\kappa = 
(V_\kappa;\in,C \cap \kappa,(W_\alpha,\pi_{\alpha \beta} 
\ \colon \alpha \leq
\beta \in C \cap \kappa)) \prec_{\Sigma_n} (V;\in,C,(W_\alpha,\pi_{\alpha \beta} 
\ \colon \alpha \leq
\beta \in C))$$ for all $\kappa \in C'$. Pick $\kappa < \lambda \in C'$, both limit
points of $C$, with
${}^\omega V_\kappa \subset V_\kappa$ and ${}^\omega V_\lambda \subset V_\lambda$
(i.e., $cf(\kappa) > \omega$ and $cf(\lambda) > \omega$).

Let ${\vec A} = (A_k \colon k<\omega)$ with $\forall k \ \exists m = m(k) \ 
A_k \in {\cal P}([\kappa]^m) \cap W$. Let $\alpha < \kappa$ be such that $A_k
\in ran(\pi_{\alpha \kappa})$ for all $k<\omega$. Set ${\bar A}_k = \pi_{\alpha
\kappa}^{-1}(A_k)$. Notice that $A_k$ is definable over
$N^\kappa$ by $$u \in A_k \Leftrightarrow {\rm \ for \ all \ but \ boundedly \
many \ } \beta \in (C \cap \kappa) \setminus
\alpha {\rm , \ } u \in \pi_{\alpha \beta}({\bar A}_k).$$
Define $B_k$ over $N^\lambda$ by 
$$u \in B_k \Leftrightarrow {\rm \ for \ all \ but \ boundedly \
many \ } \beta \in (C \cap \lambda) \setminus
\alpha {\rm , \ } u \in \pi_{\alpha \beta}({\bar A}_k).$$
Then obviously $\pi_{\kappa \lambda}(A_k) = B_k$, for all $k<\omega$. 
It is also easy to verify that $$(V_\kappa;\in,{\vec A}) \prec_{r\Sigma_3}
(V_\lambda;\in,{\vec B}).$$ (Notice that if a formula 
is $\Sigma_0(\Sigma_p)$ then it is
equivalent to a $\Sigma_p$ formula over models of $\Sigma_p$-replacement.) 

Now let $F$ be the extender derived from $\pi_{\kappa \lambda}$, and let $\nu$ be its
natural length. By our smallness assumption, $\nu < \lambda$. Let $\gamma = \nu^{+W}
< \lambda$. A straightforward
induction as in the proof of \cite{FSIT} Lemma 11.4 shows that $$({\cal
J}^W_\gamma,{\tilde F})$$ satisfies the initial segment condition, and is hence 
a premouse. But we have shown that $F$ is countably certified. Thus $F = E^W_\gamma$,
contradicting the fact that $\gamma$ is a cardinal of $W$.

\bigskip
\hfill $\square$ (\ref{Kc-is-universal})

\bigskip
Notice that \ref{main-thm} is now an immediate corollary of \ref{Kc-is-universal}
together with what we
showed in \S 1. By well-known arguments, we might in fact replace ``there's no inner
model with a Woodin cardinal'' by ``every premouse is tame,'' say, in the statement of
\ref{main-thm}.

\section{$\omega$-completeness and countable certifiability.}

We now want to discuss the relation between being $\omega$-closed and being countably
certified (in our new sense).

\begin{defn}\label{strongly-omega-closed}
Let ${\cal M}$ be an active premouse, $F$ the extender coded by ${\dot F}^{\cal M}$,
$\kappa = c.p.(F)$, and $\nu = {\dot \nu}^{\cal M} =$ the natural length of $F$. We
say that $F$ is {\em strongly $\omega$-closed} iff $\forall \ 
(a_n,X_n \colon n<\omega)$ with $$\forall n<\omega \ \exists k<\omega \ ( \ a_n \in
[\nu]^k \wedge A_n \in {\cal P}([\kappa]^k) \cap {\cal M} \ )$$
there is some transitive $N$ with $${}^\omega N \subset N \wedge V_{\nu+1} \subset N$$
such that for all $$\pi \ \colon \ {\cal M} \rightarrow_{\Sigma_2} N$$ with ${\cal M}$
countable and transitive there is 
$$\pi' \ \colon \ {\cal M} \rightarrow_{\Sigma_2} V_\kappa$$
such that $$\pi' \circ \pi^{-1} \upharpoonright \bigcup_{n<\omega} a_n \rightarrow
\kappa$$ witnesses that $F$ is $\omega$-complete w.r.t $(a_n,X_n \colon n<\omega)$,
i.e., $$\forall n<\omega \ ( \ X_n \in F_{a_n} \Rightarrow \pi' \circ \pi^{-1}(a_n)
\in X_n \ ).$$ 
\end{defn}

Recall that such $F$ is $\omega$-complete iff for all $(a_n,X_n \colon n<\omega)$ as
in \ref{strongly-omega-closed} there is an order-preserving $\tau \colon \bigcup_n a_n
\rightarrow \kappa$ with $\forall n<\omega \ ( \ X_n \in F_{a_n} \Rightarrow 
\tau(a_n) \in 
X_n \ )$. Trivially, if $F$ is strongly $\omega$-closed then $F$ is $\omega$-closed. 
Strong $\omega$-closedness requires that $\tau$ is realized as the restriction of some
$\pi' \circ \pi^{-1}$ as above. We also have the following facts, which are easy to
verify.

If $F$ is countable certified in the sense of \cite{CMIP} Defn. 1.2, then $F$ is
countably certified in the sense of \ref{certificate}, and then $F$ is strongly
$\omega$-closed. We can still run the iterability proof for countable submodels of
${\frak C}_k({\cal N}_\theta)$ if we relax the requirement that new extenders be
countably certified to that they be strongly $\omega$-complete. Of course, the new
$K^c$ is then still universal.


\begin{thebibliography}{99}
\bibitem{NFS} Ronald Jensen, {\em Addendum to A new fine structure for higher core
models}, handwritten.
\bibitem{FSIT} Bill Mitchell and John Steel, {\em Fine structure and iteration trees},
LNL $\sharp 3$.
\bibitem{K} Ralf Schindler, {\em The core model for almost linear iterations},
submitted.
\bibitem{CMIP} John Steel, {\em The core model iterability problem}, LNL $\sharp 8$.
\end{thebibliography}
\end{document}